\documentclass[preprint, 3p, numbers]{elsarticle}

\usepackage{lineno,hyperref}

\usepackage{amsmath}
\usepackage{amssymb}
\usepackage{amsthm}

\usepackage{mathptmx}      

\usepackage{natbib}
\usepackage{bibentry} 

\usepackage{subfig}
\usepackage{graphicx}

\usepackage[active]{srcltx} 

\usepackage[bbgreekl]{mathbbol}
\usepackage{bm} 
\usepackage{MnSymbol} 

\usepackage{units}
\usepackage{tensor}
\usepackage{accents}

\DeclareMathOperator{\Tr}{Tr}

\newcommand{\bigo}[1]{\ensuremath{O\left(#1 \right)}}

\newcommand{\exponential}[1]{\ensuremath{{\mathrm e}^{#1}}}

\newcommand{\iunit}{\ensuremath{\mathrm{i}}}

\newcommand{\bydefinition}{\mathrm{def}}

\newcommand{\dimless}[1]{#1^\star}

\newcommand{\diff}{\mathrm{d}}

\makeatletter
\@ifpackageloaded{bm}%
{%
}{%
\relax
}
\makeatother
\newcommand{\tensorq}[1]{\ensuremath{\mathbb{#1}}}      

\makeatletter
\@ifpackageloaded{bm}%
{%
}{%

}

\@ifpackageloaded{bm}%
{%
 
}{%

}

\@ifpackageloaded{bm}%
{%
}{%

}
\makeatother

\newcommand{\diracdelta}{\delta}

\newcommand{\R}{\ensuremath{{\mathbb R}}}

\newcommand{\N}{\ensuremath{{\mathbb N}}}
\newcommand{\Z}{\ensuremath{{\mathbb Z}}}

\newcommand{\dd}[2]{\ensuremath{\frac{\diff {#1}}{\diff {#2}}}}

\newcommand{\ddd}[2]{\ensuremath{\frac{\diff^2 {#1}}{\diff {#2}^2}}}

\makeatletter
\@ifpackageloaded{tensor}
{

}{%

}
\makeatother

\makeatletter
\@ifpackageloaded{tensor}
{

}{%

}
\makeatother

\newcommand{\absnorm}[1]{\ensuremath{\left|#1\right|}}

\makeatletter
\@ifundefined{volume}{%
}%
{%
}
\makeatother

\makeatletter
\@ifpackageloaded{MnSymbol} 
{
 
}{%
 
}
\makeatother

\numberwithin{equation}{section}
\let\cite\citet

\begin{document}

\begin{frontmatter}

\title{A note on parametric resonance induced by a singular parameter modulation}

\author[charles]{Dalibor Pra\v{z}\'ak}
\ead{prazak@karlin.mff.cuni.cz}

\author[charles]{V\'{\i}t Pr\r{u}\v{s}a\fnref{myfootnote1}\corref{mycorrespondingauthor}}
\fntext[myfootnote1]{V\'{\i}t Pr\r{u}\v{s}a and Dalibor Pra\v{z}\'ak thank the Czech Science Foundation, grant number~20-11027X, for its support.}
\ead{prusv@karlin.mff.cuni.cz}

\author[charles]{Karel T\r{u}ma\fnref{myfootnote2}}
\ead{ktuma@karlin.mff.cuni.cz}
\fntext[myfootnote2]{Karel T\r{u}ma was supported by Charles University research program No. UNCE/SCI/023.}

\address[charles]{Charles University, Faculty of Mathematics and Physics\\
Sokolovsk\'a 83, Praha, CZ 186 75, Czech Republic}

\begin{abstract}
  We investigate the classical problem of motion of a mathematical pendulum with an oscillating pivot. This simple mechanical setting is frequently used as the prime example of a system exhibiting the parametric resonance phenomenon, which manifests itself by surprising stabilisation/destabilisation effects. In the classical case the pivot oscillations are described by a cosine wave, and the corresponding stability analysis requires one to investigate the behaviour of solutions to the Mathieu equation. This is not a straightforward procedure, and it does not lead to exact and simple analytical results expressed in terms of elementary functions. Consequently, the explanation of the parametric resonance phenomenon can be in this case obscured by the relatively involved technical calculations. We show that the stability analysis is much easier if one considers the pivot motion described by a non-smooth function---a triangular or a nearly rectangular wave. The non-smooth pivot motion leads to the presence of singularities (Dirac distributions) in the corresponding Mathieu type equation, which seemingly further complicates the analysis. Fortunately, this is only a minor technical difficulty. Once the mathematical setting for the non-smooth forcing is settled down, the corresponding stability diagram is indeed straightforward to obtain, and the stability boundaries are, unlike in the classical case, given in terms of simple analytical formulae involving only elementary functions.


\end{abstract}


\end{frontmatter}


\section{Introduction}
\label{sec:introduction}
Parametric resonance is a well-known dynamical phenomenon that can take place in systems with periodically varying parameters. The core observation is that the modulation of system parameters---rather than a direct external forcing as in the standard resonance phenomenon---can have striking stabilising/destabilising effect on the dynamical behaviour of the system.

The classical manifestation of the parametric resonance is the stabilisation of the inverted pendulum by the means of an oscillating pivot, which is the problem first analysed by~\cite{stephenson.a:on*2} and later independently by~\cite{kapitza.p:dynamic}; see~\cite{acheson.d:pendulum} for further historical remarks. The fact that the pendulum can be stabilised in the upright position by fast vertical pivot oscillations is from the naive point of view surprising, and it is frequently discussed in introductory texts on mechanical vibrations, see~\cite{hartog.jp:mechanical}, \cite{landau.ld.lifshitz.em:course*3} or~\cite{nayfeh.ah:introduction,nayfeh.ah:perturbation} to name a few, and on theory of ordinary differential equations, see, for example, \cite{jordan.dw.smith.p:nonlinear}.

The explanation of the phenomenon is usually given in terms of the analysis of Mathieu equation
\begin{equation}
  \label{eq:39}
  \ddd{\dimless{\theta}}{\dimless{t}}
  +
  \left(
    \alpha
    +
    \beta
    \cos \dimless{t}
  \right)
  \dimless{\theta}
  =
  0
  ,
\end{equation}
where $\alpha$ and $\beta$ are some parameters, and $\dimless{\theta}$ describes the angular displacement of the pendulum, see below for details. The behaviour of solutions to~\eqref{eq:39} with respect to parameters $\alpha$ and $\beta$ is typically analysed using the standard Floquet theory followed by a perturbation technique, that is employed either directly on the level of the equation itself, see, for example, \cite{nayfeh.ah:introduction,nayfeh.ah:perturbation}, or indirectly in the manipulations with the so-called Hill determinant, see~\cite{hill.gw:on} and modern discussion in~\cite{morse.pm.feshbach.h:methods}, \cite{phelps.fm.hunter.jh:analytical} or~\cite{jordan.dw.smith.p:nonlinear} to name a few. The other possibility is to employ a heuristics based on the slow/fast time scale separation and the construction of an effective potential, which is the method originally applied by~\cite{kapitza.p:dynamic}; for further references in this regard see newer works such as~\cite{butikov.ei:on} or \cite{artstein.z:pendulum} and references therein. These methods either give analytical results valid only in a small amplitude and/or high frequency approximation  (small parameter~$\beta$, small parameter~$\alpha$), or require one to elaborately manipulate the Hill determinant, while the evaluation of the determinant anyway finally resorts to yet another approximation procedure, see~\cite{jordan.dw.smith.p:nonlinear}.

Consequently, the analysis of the parametric resonance phenomenon for inverted pendulum remains either only on a heuristic level or it is obscured by relatively involved technical calculations. In principle, this is not a problem since the heuristics suffice to give an insight into the phenomenon, and the rigorous analysis of Mathieu equation---as well as the more general Hill equation---is nowadays a well-developed field, see~\cite{mclachlan.nw:theory} and~\cite{magnus.wws:hills}. However, since the parametric resonance phenomenon is an important phenomenon in technical practice, see~\cite{champneys.a:dynamics} and references therein, and since it is a matter of everlasting curiosity, see, for example, \cite{keller.jb:ponytail}, there is a relentless effort to provide a simple, intuitive and rigorous enough explanation of the phenomenon without the need to resort to cumbersome technical calculations.

Such attempts are, for example, based on a slight reformulation of the problem, see~\cite{hartog.jp:mechanical} or~\cite{levi.m.weckesser.w:stabilization}, or on various physical arguments see, for example, a more recent discussion in~\cite{butikov.ei:on,butikov.ei:improved,butikov.ei:analytical} and~\cite{artstein.z:pendulum}. In what follows we provide yet another instructive reformulation of the parametric resonance problem for the pendulum with an oscillating pivot, while the presented reformulation~\emph{allows one to explicitly carry out all the necessary calculations in terms of simple expressions involving only elementary functions}.

The key modification leading to a simple subsequent analysis is the choice of a suitable pivot motion modulation. It turns out that the \emph{triangular wave} instead of the classical \emph{cosine wave} is the suitable choice. (Note that the original analysis by~\cite{stephenson.a:on*2} partially went in this direction as well.) We emphasise that we deal with a triangular wave for the \emph{pivot motion}, we are \emph{not} dealing with a triangular wave in the \emph{time-periodic coefficient} in the Mathieu equation. The latter case is frequently studied in the literature and it is of no interest here.

The fact that we work with a \emph{non-smooth pivot motion} however leads to a technical difficulty. The reason is that the \emph{acceleration} of the pivot motion is given in terms of generalised functions. In particular, the classical Mathieu equation~\eqref{eq:39} is replaced by
\begin{equation}
  \label{eq:40}
    \ddd{\dimless{\theta}}{\dimless{t}}
  +
  \left(
    \alpha
    +
    \beta
    \left(
      \diracdelta_{\dimless{t} - \left(2n+1\right) \pi}
      -
      \diracdelta_{\dimless{t} - 2n \pi}
    \right)
  \right)
  \dimless{\theta}
  =
  0
  ,
\end{equation}
where $\diracdelta_{\dimless{t} - a}$ denotes the Dirac distribution with respect to the variable $\dimless{t}$ located at point $a$. (See Section~\ref{sec:governing-equations} for details and detailed discussion of the notation.) This means that we face the problem of multiplication of the solution $\dimless{\theta}$ with the generalised function $\diracdelta_{\dimless{t} - a}$, which is an operation not defined in the classical theory of distributions, see~\cite{schwartz.l:sur}.  This technical difficulty is in the present case easy to overcome, and the remaining calculations necessary for the stability analysis are straightforward.

After the analysis of the resonance induced by the triangular wave, we proceed with yet another non-smooth pivot motion, namely with a (nearly) \emph{rectangular wave}. This is a very interesting case of non-smooth pivot motion, since it allows one to challenge the conventional wisdom, namely the fact that ``practicing engineers often think of this \emph{subharmonic instability} close to $\alpha = \frac{1}{4}$ [\dots] as being \emph{the} halmark of parametric resonance'', see~\cite{champneys.a:dynamics}.

\section{Governing equations}
\label{sec:governing-equations}

\begin{figure}[h]
\centering
\includegraphics[width=0.25\textwidth]{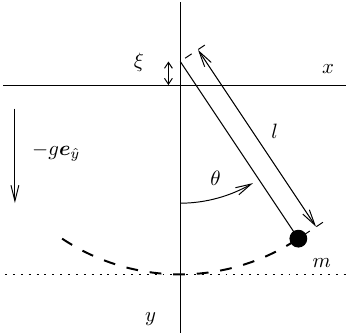}
\caption{Mathematical pendulum with an oscillating pivot.}  
\label{fig:inverted-pendulum-intro}
\end{figure}

Using the standard techniques of analytical mechanics, see, for example, \cite{meirovitch.l:fundamentals}, the governing equations for pendulum of length $l$ swinging in a homogeneous gravitational field with the gravitational acceleration $g$ and with the pivot localised at the vertical position $\xi(t)$ are found to be
\begin{equation}
  \label{eq:1}
  ml^2 \ddd{\theta}{t}
  +
  ml
  \left(
    \ddd{\xi}{t}
    +
    g
  \right)
  \sin \theta
  =
  0
  ,
\end{equation}
where $\theta$ denotes the angle between the vertical axis and the pendulum, see Figure~\ref{fig:inverted-pendulum-intro}. If we assume that the pivot is oscillating with a minimal period $T$ and the corresponding angular frequency $\Omega = \frac{2 \pi}{T}$, and if we introduce the dimensionless time $\dimless{t}$ as $\dimless{t}=_{\bydefinition} \Omega t$, that is $\dimless{t}=_{\bydefinition} 2 \pi \frac{t}{T}$, then the governing equation~\eqref{eq:1} reduces to
\begin{equation}
  \label{eq:6}
  \ddd{\dimless{\theta}}{\dimless{t}}
  +
  \left(
    \frac{1}{l}
    \ddd{\xi}{\dimless{t}}
    +
    \alpha
  \right)
  \sin \dimless{\theta}
  =
  0
  ,
\end{equation}  
where $\alpha =_{\bydefinition} \frac{\omega_0^2}{\Omega^2}$ and $\omega_0 = _{\bydefinition} \sqrt{\frac{g}{l}}$ denotes the natural frequency of the pendulum.

\begin{figure}[h]
  \centering
  \subfloat[Triangular wave.\label{fig:pivot-motion-a}]{\includegraphics[width=0.32\textwidth]{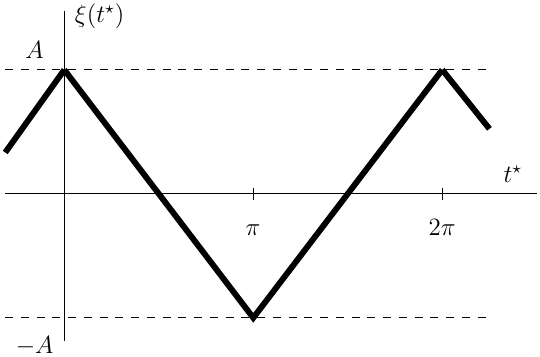}}
  \qquad
  \subfloat[Approximated rectangular wave.\label{fig:pivot-motion-b}]{\includegraphics[width=0.32\textwidth]{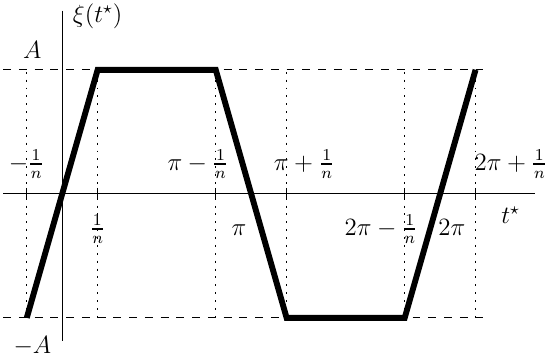}}
  \caption{Pivot motion~$\xi(\dimless{t})$, dimensionless time.}
  \label{fig:pivot-motion}
\end{figure}

We consider the response of the system to two nonstandard pivot motions $\xi(\dimless{t})$: the pivot motion described by a \emph{triangular wave}, see Figure~\ref{fig:pivot-motion-a}, and the pivot motion described by an approximation of a \emph{rectangular wave}, see Figure~\ref{fig:pivot-motion-b}. Both nonstandard pivot motions are in principle realisable by a simple mechanical device. For example, running the appropriately configured old mechanical beauty---the \emph{harmonic analyser}, see~\cite{michelson.aa.stratton.sw:new}---in the reverse direction will in principle lead to the desired pivot motion. The outcomes of analysis for both nonstandard pivot motions are compared to the classical setting wherein the pivot motion is described by a~\emph{cosine wave}.

The behaviour of solutions to~\eqref{eq:6} is typically investigated in the \emph{linearised setting}. The term $\sin \dimless{\theta}$ is linearised as $\sin \dimless{\theta} \approx \dimless{\theta}$ (pendulum is swinging close to the pendent position) or as $\sin \dimless{\theta} \approx - \dimless{\zeta}$ (pendulum is swinging close to the upright position $\dimless{\theta} = \dimless{\zeta} + \pi$, $\dimless{\zeta} \approx 0$). We stick to the linearised setting as well. Regarding the nonlinear setting and the standard cosine wave case we refer the interested reader to~\cite{kidachi.h.onogi.h:note} and references therein. Note that the upright position case is formally identical to the pendent position case, the only difference is that in the latter case we work with \emph{negative} values of parameters~$\alpha$ and~$\beta$.

\begin{figure}[h]
  \centering
  \subfloat[Triangular wave.\label{fig:pivot-acceleration-a}]{\includegraphics[width=0.32\textwidth]{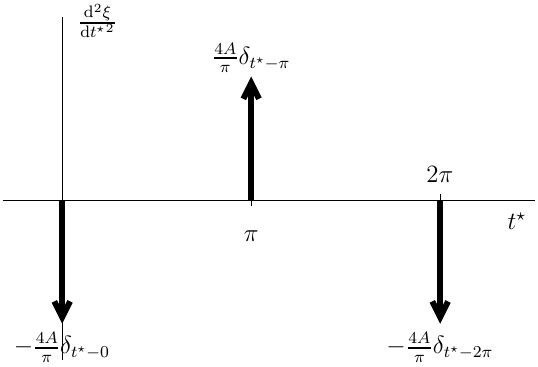}}
  \qquad
  \subfloat[Approximated rectangular wave.\label{fig:pivot-acceleration-b}]{\includegraphics[width=0.32\textwidth]{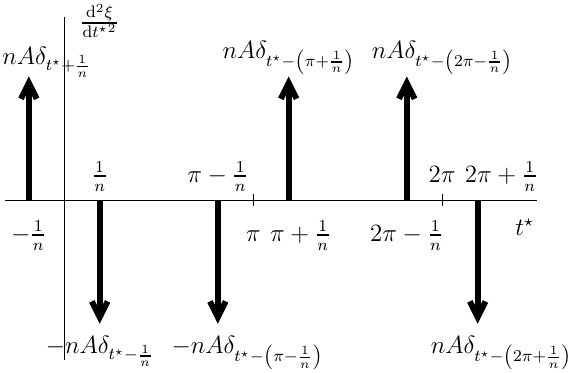}}
  \caption{Pivot acceleration~$\ddd{\xi}{\dimless{t}}$, dimensionless time.}
  \label{fig:pivot-acceleration}
\end{figure}

In the standard setting one uses the cosine wave for the pivot motion, $\xi =_{\bydefinition} A \cos \left(\Omega t \right)$, and the linearisation leads to the Mathieu equation~\eqref{eq:39} for $\dimless{\theta}$ and $\dimless{\zeta}$ respectively. In our generalised setting we consider the pivot motion~$\xi$ given by a triangular wave~\eqref{eq:2} or by an approximated rectangular wave~\eqref{eq:44}, but in the end we still work with a linear ordinary differential equation with time periodic coefficients. This means that we can still in principle use the stability analysis based on the standard Floquet theory, see~\cite{floquet.g:sur} and a modern discussion thereof in, for example, \cite{jordan.dw.smith.p:nonlinear}. However, we need to make necessary amendments of Floquet theory to the setting of generalised functions. Indeed, since the pivot motion $\xi$ is not smooth, the second derivative $\ddd{\xi}{\dimless{t}}$ in the governing equation~\eqref{eq:6} involves singular terms, namely the Dirac distributions centered at given points, see Figure~\ref{fig:pivot-motion} and Figure~\ref{fig:pivot-acceleration} respectively. This is however a minor technical difficulty, and we discuss it below.

In Floquet theory the stability analysis essentially boils down to the evaluation of the trace of the monodromy matrix~$\tensorq{E}$, which is the fundamental matrix $\tensorq{\Phi}$ evaluated at time $\dimless{t} = 2 \pi$, see the detailed discussion in~\cite{jordan.dw.smith.p:nonlinear}. The monodromy matrix $\tensorq{E}$ for the system forced by a triangular/rectangular wave is simply assembled using the matrix product of fundamental matrices for a single jump forcing and the free (unforced) regime. Consequently, before we proceed with full analysis of triangular/rectangular wave, it is worthwhile to study the response of the system to a single jump forcing and the response of the system in the free regime.

\section{System response in simple cases}
\label{sec:syst-resp-simple}

\subsection{System response to a single jump}
\label{sec:syst-resp-single}
We first investigate what happens if the system governed by the equation
\begin{equation}
  \label{eq:11}
    \ddd{\dimless{\theta}}{\dimless{t}}
  +
  \left(
    \alpha
    +
    \Gamma
    \diracdelta_{\dimless{t} - \dimless{t}_{\text{jump}}}
  \right)
  \dimless{\theta}
  =
  0
  ,
\end{equation}
where $\Gamma \in \R$ is a fixed number, crosses the singularity at $\dimless{t}=\dimless{t}_{\text{jump}}$. It is straightforward to verify that traversing the singularity at $\dimless{t}=\dimless{t}_{\text{jump}}$ induces a jump in the derivative of the function $\dimless{\theta}$, while the function $\dimless{\theta}$ itself remains continuous across the jump, that is
\begin{equation}
  \label{eq:12}
    \lim_{\dimless{t} \to \dimless{t}_{\text{jump}}} \dimless{\theta}=\lim_{\dimless{t} \to \dimless{t}_{\text{jump}} -} \dimless{\theta}
    =
      \lim_{\dimless{t} \to \dimless{t}_{\text{jump}} +} \dimless{\theta}
      ,
    \qquad
    \lim_{\dimless{t} \to \dimless{t}_{\text{jump}} +} \dd{\dimless{\theta}}{\dimless{t}}
    =
      \lim_{\dimless{t} \to \dimless{t}_{\text{jump}} -} \dd{\dimless{\theta}}{\dimless{t}}
      -
      \Gamma \left. \dimless{\theta} \right|_{\dimless{t} = \dimless{t}_{\text{jump}}}.
\end{equation}
If we rewrite conditions~\eqref{eq:12} in a matrix form we see that
\begin{equation}
  \label{eq:14}
  \left.
  \begin{bmatrix}
    \dimless{\theta}  \\
    \dd{\dimless{\theta}}{\dimless{t}}
  \end{bmatrix}
  \right|_{\dimless{t} \to \dimless{t}_{\text{jump}} +}
  =
  \begin{bmatrix}
    1 & 0 \\
    -\Gamma & 1
  \end{bmatrix}
  \left.
    \begin{bmatrix}
      \dimless{\theta}  \\
      \dd{\dimless{\theta}}{\dimless{t}}
    \end{bmatrix}
  \right|_{\dimless{t} \to \dimless{t}_{\text{jump}} -}
  ,
\end{equation}
hence the fundamental matrix $\tensorq{\Phi}_{\dimless{t}_{\text{jump}}- \to \dimless{t}_{\text{jump}}+}$ that transfers the system across the singularity reads 
\begin{equation}
  \label{eq:transfer-matrix-jump}
  \tensorq{\Phi}_{\dimless{t}_{\text{jump}}- \to \dimless{t}_{\text{jump}}+}
  =
  \begin{bmatrix}
    1 & 0 \\
    -\Gamma & 1
  \end{bmatrix}
  .
\end{equation}

Note however that our argument is rather sloppy from the rigorous point of view. In equation~\eqref{eq:11} we are multiplying two generalised functions/distributions, namely the sought function $\dimless{\theta}$ and the Dirac distribution $\diracdelta_{\dimless{t} - \dimless{t}_0}$, which is not possible in the classical theory of distributions, see~\cite{schwartz.l:sur}. This theory of generalised functions allows one to multiply a distribution \emph{only} by a smooth function. In other cases the multiplication operation is not defined.

Consequently, we need to employ another viewpoint regarding the generalised functions/distributions. We need a structure that allows us to transparently and simultaneously handle discontinuity, differentiation and nonlinearity. This can be done for example in the so-called Colombeau algebra, see especially~\cite{colombeau.j:new,colombeau.j:elementary,colombeau.j:multiplication}, \cite{biagioni.ha:nonlinear}, \cite{rosinger.ee:generalized} or \cite{grosser.m.farkas.e.ea:on}; for applications of Colombeau algebra in physics see, for example, \cite{colombeau.jf.le-roux.ay:multiplications}, \cite{grosser.m.kunzinger.m.ea:geometric}, \cite{steinbauer.r.vickers.ja:use}, \cite{aragona.j.colombeau.jf.ea:nonlinear}, \cite{todorov.td:steady-state}, \cite{rehor.m.pusa.v.ea:on} and \cite{prusa.v.rajagopal.kr:on,prusa.v.rehor.m.ea:colombeau*1} to name a few. We shall however not go into the technical details, we will be content with the claim that in the present case all manipulations work as expected, meaning that the product of a continuous function and the Dirac distribution yields the value of the function at the corresponding point, and that all the required manipulations can be---if necessary---formalised.  

\subsection{System response in the free regime}
\label{sec:system-response-free}
Next we investigate what happens if the system governed by the equation
$
\ddd{\dimless{\theta}}{\dimless{t}}
  +
  \alpha
  \dimless{\theta}
  =
  0
$
  ,
evolves from time $\dimless{t}_0$ to time $\dimless{t}_0+ \dimless{\tau}$. Elementary calculation yields that 
\begin{equation}
  \label{eq:43}
  \left.
  \begin{bmatrix}
    \dimless{\theta}  \\
    \dd{\dimless{\theta}}{\dimless{t}}
  \end{bmatrix}
  \right|_{\dimless{t} = \dimless{t}_0 + \dimless{\tau}}
  =
  \left.
    \left(
      \exponential
      {
        \begin{bmatrix}
          0 & 1 \\
          -\alpha & 0
        \end{bmatrix}
        (\dimless{t} - \dimless{t}_0)
      }
    \right)
  \right|_{\dimless{t} = \dimless{t}_0 + \dimless{\tau}}
  \left.
    \begin{bmatrix}
      \dimless{\theta}  \\
      \dd{\dimless{\theta}}{\dimless{t}}
    \end{bmatrix}
  \right|_{\dimless{t} = \dimless{t}_0}
  =
    \begin{bmatrix}
      \cos \left(  \sqrt{\alpha} \dimless{\tau}\right) &   \frac{1}{\sqrt{\alpha}} \sin \left(  \sqrt{\alpha} \dimless{\tau}\right) \\
      -\sqrt{\alpha} \sin \left(  \sqrt{\alpha} \dimless{\tau}\right) & \cos \left(  \sqrt{\alpha} \dimless{\tau}\right)
    \end{bmatrix}
  \left.
    \begin{bmatrix}
      \dimless{\theta}  \\
      \dd{\dimless{\theta}}{\dimless{t}}
    \end{bmatrix}
  \right|_{\dimless{t} = \dimless{t}_0}
  .
\end{equation}
Consequently the fundamental matrix $\tensorq{\Phi}_{\dimless{t}_0 \to \dimless{t}_0 + \dimless{\tau}}$ that transfers the system in the free regime from time $\dimless{t}_0$ to time $ \dimless{t}_0 + \dimless{\tau}$  reads 
\begin{equation}
  \label{eq:transfer-matrix-free}
  \tensorq{\Phi}_{\dimless{t}_0 \to  \dimless{t}_0 + \dimless{\tau}}
  =
  \begin{bmatrix}
    \cos \left(  \sqrt{\alpha} \dimless{\tau}\right) &   \frac{1}{\sqrt{\alpha}} \sin \left(  \sqrt{\alpha} \dimless{\tau}\right) \\
    -\sqrt{\alpha} \sin \left(  \sqrt{\alpha} \dimless{\tau}\right) & \cos \left(  \sqrt{\alpha} \dimless{\tau}\right)
  \end{bmatrix}
  .
\end{equation}

\section{Pivot motion given as a triangular wave}
\label{sec:pivot-motion-given}
Having obtained the fundamental matrix for the free regime over the time interval of length $\tau$ and the fundamental matrix for the singularity crossing, we can proceed with the identification of the monodromy matrix for the pivot motion described by a triangular wave.

\subsection{Triangular wave -- governing equations}
\label{sec:triang-wave-govern}
A triangular wave with amplitude $A$ and minimal period $2\pi$ in the dimensionless time, see Figure~\ref{fig:pivot-motion-a}, is given by the function $\xi$ specified by the formula
\begin{equation}
  \label{eq:2}
  \xi(\dimless{t})
  =
  \begin{cases}
    A - \frac{2A}{\pi} \dimless{t}, & \dimless{t} \in \left[0, \pi \right], \\
    -3A + \frac{2A}{\pi} \dimless{t}, & \dimless{t} \in \left(\pi, 2\pi \right).
  \end{cases}
\end{equation}
The first and second derivatives of $\xi$ are then---informally---given by the formulae
\begin{subequations}
  \label{eq:4}
  \begin{align}
  \label{eq:3}
    \dd{\xi}{\dimless{t}}
    &=
      \frac{2A}{\pi}
      \left(
      -
      \chi_{\left[0, \pi\right]}
      +
      \chi_{\left(\pi, 2\pi\right)}
      \right)
      , \\
    \label{eq:41}
    \ddd{\xi}{\dimless{t}}
    &=
      \frac{4A}{\pi}
      \left(
      \diracdelta_{\dimless{t} - \left(2l+1\right) \pi}
      -
      \diracdelta_{\dimless{t}- 2l\pi}
      \right)
      ,
  \end{align}
\end{subequations}
where $l \in \Z$, and where the symbol $\diracdelta_{\dimless{t}-\dimless{t}_0}$ denotes the Dirac distribution centered at point $\dimless{t}=\dimless{t}_0$, and the symbol $\chi_{\left[0, \pi \right]}$ denotes the characteristic function of the corresponding interval. The acceleration $\ddd{\xi}{\dimless{t}}$ is schematically shown in Figure~\ref{fig:pivot-acceleration-a}, and it consists of positive/negative impulses located at points $0$ and $\pi$ with period $2\pi$. We also note that impulsive coefficients are sometimes treated in the literature, see, for example, \cite{richards.ja:analysis}, who deals with a sequence of unidirectional impulses. In our setting it is important that the system is subject to the sequence of \emph{alternating} impulses. This case seems to be not dealt with in the literature so far.

Substituting~\eqref{eq:41} into the governing equation~\eqref{eq:6} yields
\begin{equation}
  \label{eq:7}
  \ddd{\dimless{\theta}}{\dimless{t}}
  +
  \left(
    \alpha
    +
    \beta
    \left(
      \diracdelta_{\dimless{t} - \left(2l+1\right) \pi}
      -
      \diracdelta_{\dimless{t} - 2l\pi}
    \right)
  \right)
  \sin \dimless{\theta}
  =
  0
  ,
\end{equation}
where the dimensionless parameter $\beta$ is given by the formula $\beta =_{\bydefinition} \frac{4A}{\pi l}$. The linearisation of~\eqref{eq:7} with respect to $\dimless{\theta}$ then leads to
\begin{subequations}
  \label{eq:8}
  \begin{align}
    \label{eq:9}
  \ddd{\dimless{\theta}}{\dimless{t}}
  +
  \left(
    \alpha
    +
    \beta
    \left(
      \diracdelta_{\dimless{t} - \left(2l+1\right) \pi}
      -
      \diracdelta_{\dimless{t} - 2l\pi}
    \right)
  \right)
  \dimless{\theta}
  &=
  0
  ,
    \\
    \label{eq:10}
  \ddd{\dimless{\theta}}{\dimless{t}}
  -
  \left(
    \alpha
    +
    \beta
    \left(
      \diracdelta_{\dimless{t} - \left(2l+1\right) \pi}
      -
      \diracdelta_{\dimless{t} - 2l\pi}
    \right)
  \right)
  \dimless{\theta}
  &=
  0
  ,
  \end{align}
\end{subequations}
depending whether we linearise in the vicinity of the point $0$ or $\pi$. In the former case---the pendulum is hanging down---we get~\eqref{eq:9}, while in the latter case---the pendulum is in the inverted position---we get~\eqref{eq:10}. This formally allows us to investigate positive and negative values of parameters $\alpha$ and $\beta$. Positive values then correspond to the pendulum in the pendent position, while negative values correspond to the pendulum in the upright position. 

\subsection{Monodromy matrix}
\label{sec:monodromy-matrix}
The monodromy matrix $\tensorq{E}$ for the system of interest~\eqref{eq:9} is the fundamental matrix evaluated over the whole period $\tensorq{\Phi}_{0- \to 2 \pi-}$, and it is given by the product of the fundamental matrices for the free regime and the fundamental matrices for the single jump of height $\pm\beta$ located at time instants $\dimless{t}_{\text{jump}}=0$ and $\dimless{t}_{\text{jump}}=\pi$, see~\eqref{eq:transfer-matrix-free} and~\eqref{eq:transfer-matrix-jump}. We get
\begin{multline}
  \label{eq:16}
  \tensorq{E}
  =
  \tensorq{\Phi}_{0- \to 2 \pi-}
  =
  \tensorq{\Phi}_{\pi+ \to 2 \pi-}
  \tensorq{\Phi}_{\pi- \to \pi+}
  \tensorq{\Phi}_{0+ \to \pi-}
  \tensorq{\Phi}_{0- \to 0+}
  \\
  =
  \begin{bmatrix}
    \cos \left(  \sqrt{\alpha} \pi \right) &   \frac{1}{\sqrt{\alpha}} \sin \left(  \sqrt{\alpha} \pi \right) \\
    -\sqrt{\alpha} \sin \left(  \sqrt{\alpha} \pi \right) & \cos \left(  \sqrt{\alpha} \pi \right)
  \end{bmatrix}
  \begin{bmatrix}
    1 & 0 \\
    -\beta  & 1
  \end{bmatrix}
  \begin{bmatrix}
    \cos \left(  \sqrt{\alpha} \pi \right) &   \frac{1}{\sqrt{\alpha}} \sin \left(  \sqrt{\alpha} \pi \right) \\
    -\sqrt{\alpha} \sin \left(  \sqrt{\alpha} \pi \right) & \cos \left(  \sqrt{\alpha} \pi \right)
  \end{bmatrix}
  \begin{bmatrix}
    1 & 0 \\
    \beta & 1
  \end{bmatrix}
  .
\end{multline}
We see that
$
\Tr \tensorq{\Phi}_{0- \to 2 \pi-}
  =
  2 \cos \left(2 \pi \sqrt{\alpha} \right) - \frac{\beta^2}{\alpha} \sin^2 \left(\pi \sqrt{\alpha}\right)
$,
which upon application of trigonometric identities yields
\begin{equation}
  \label{eq:18}
  \Tr \tensorq{\Phi}_{0- \to 2 \pi-}
  =
  2 \cos^2 \left(\pi \sqrt{\alpha} \right) - \left( 2 + \frac{\beta^2}{\alpha}  \right)\sin^2 \left(\pi \sqrt{\alpha}\right),
\end{equation}
which is the most convenient form for the ongoing manipulations. If $\alpha<0$ we can rewrite~\eqref{eq:18} in terms of hyperbolic functions---we interpret $\sqrt{\alpha}$ as $\iunit \sqrt{\absnorm{\alpha}}$, and we use identities $\cosh x = \cos \left( \iunit x \right)$ and $\iunit \sinh x = \sin \left( \iunit x \right)$.

\subsection{Search for periodic solutions}
\label{sec:search-peri-solut}
Following the standard analysis, see~\cite{jordan.dw.smith.p:nonlinear}, we know that the system~\eqref{eq:8} possesses a periodic solution provided that $\Tr \tensorq{\Phi}_{0- \to 2 \pi-} = \pm 2$, and we know that the implicitly defined curves $\Tr \tensorq{\Phi}_{0- \to 2 \pi-} = \pm 2$ are the boundaries of stable/unstable regions in the parameter space $(\alpha, \beta) \in \R^2$. 

\subsubsection{Solution of equation $\Tr \tensorq{\Phi}_{0- \to 2 \pi-} = 2$}
\label{sec:solution-equation-tr}

If we use expression~\eqref{eq:18} for the trace of the monodromy matrix, we see that for $\alpha >0$ the equation $\Tr \tensorq{\Phi}_{0- \to 2 \pi-} = 2$ reduces to
\begin{equation}
  \label{eq:19}
  - \left( 4 + \frac{\beta^2}{\alpha} \right) \sin^2 \left(\pi \sqrt{\alpha}\right) = 0,
\end{equation}
while in the case of $\alpha < 0$ we get
\begin{equation}
  \label{eq:20}
  \left(4 - \frac{\beta^2}{\absnorm{\alpha}} \right) \sinh^2 \left(\pi \sqrt{\absnorm{\alpha}}\right) = 0.
\end{equation}

\subsubsection{Solution of equation for $\Tr \tensorq{\Phi}_{0- \to 2 \pi-} = -2$}
\label{sec:solution-equation-tr-1}
For $\alpha >0$ the equation reduces to
\begin{equation}
  \label{eq:21}
  \absnorm{\beta} = 2 \sqrt{\alpha} \absnorm{\frac{\cos \left( \pi \sqrt{\alpha} \right)}{\sin \left( \pi \sqrt{\alpha} \right)}},
\end{equation}
while in the case of $\alpha < 0$ we get
\begin{equation}
  \label{eq:22}
  \absnorm{\beta} = 2 \sqrt{\absnorm{\alpha}} \absnorm{\frac{\cosh \left( \pi \sqrt{\absnorm{\alpha}} \right)}{\sinh \left( \pi \sqrt{\absnorm{\alpha}} \right)}}.
\end{equation}

\subsubsection{Ince--Strutt diagram}
\label{sec:ince-strutt-diagram}
\begin{figure}[h]
  \centering
  \subfloat[Global view.\label{fig:triangular-wave-plots-a}]{\includegraphics[height=0.32\textwidth]{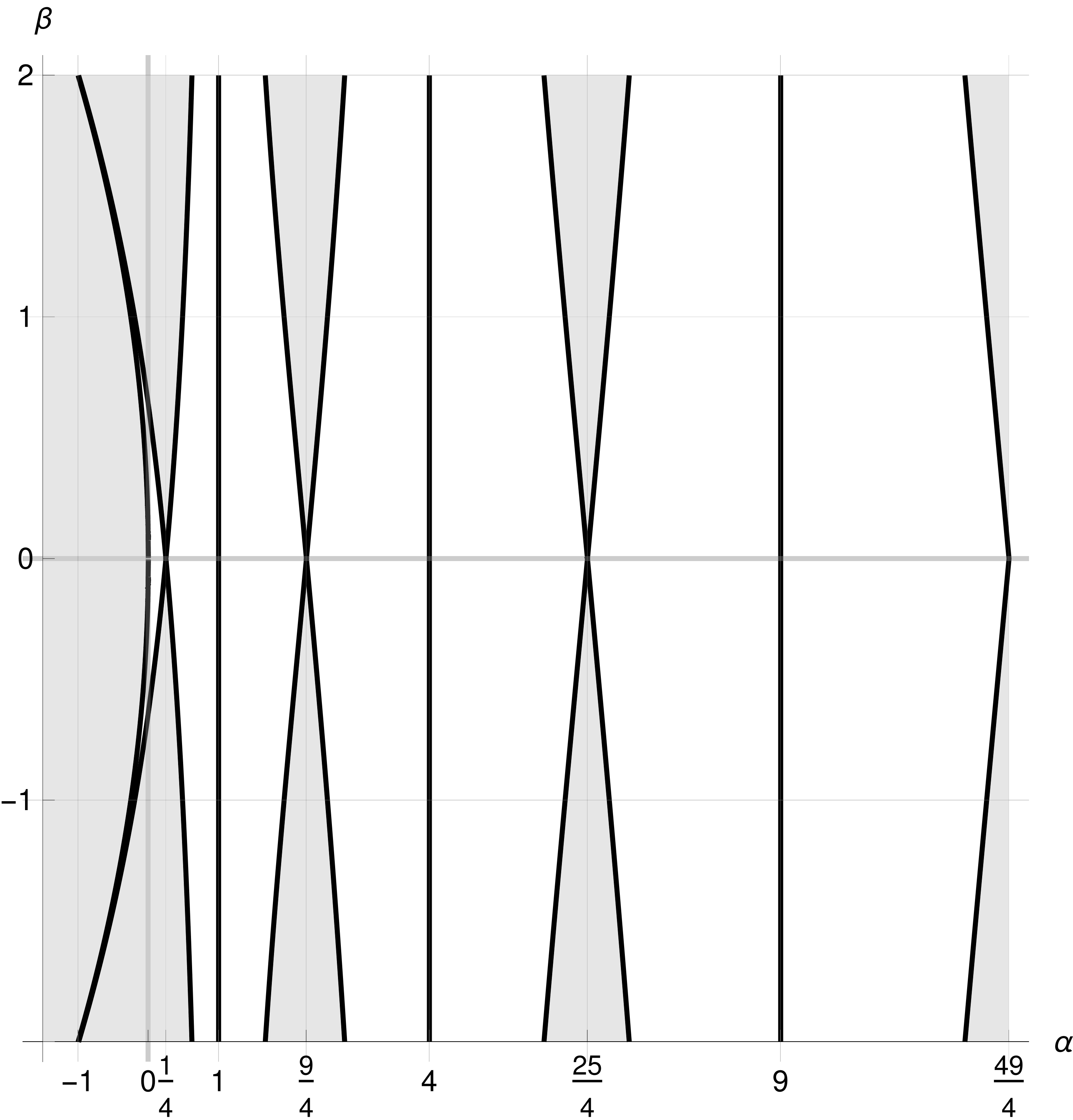}}
  \quad
  \subfloat[Detailed view.\label{fig:triangular-wave-plots-b}]{\includegraphics[height=0.32\textwidth]{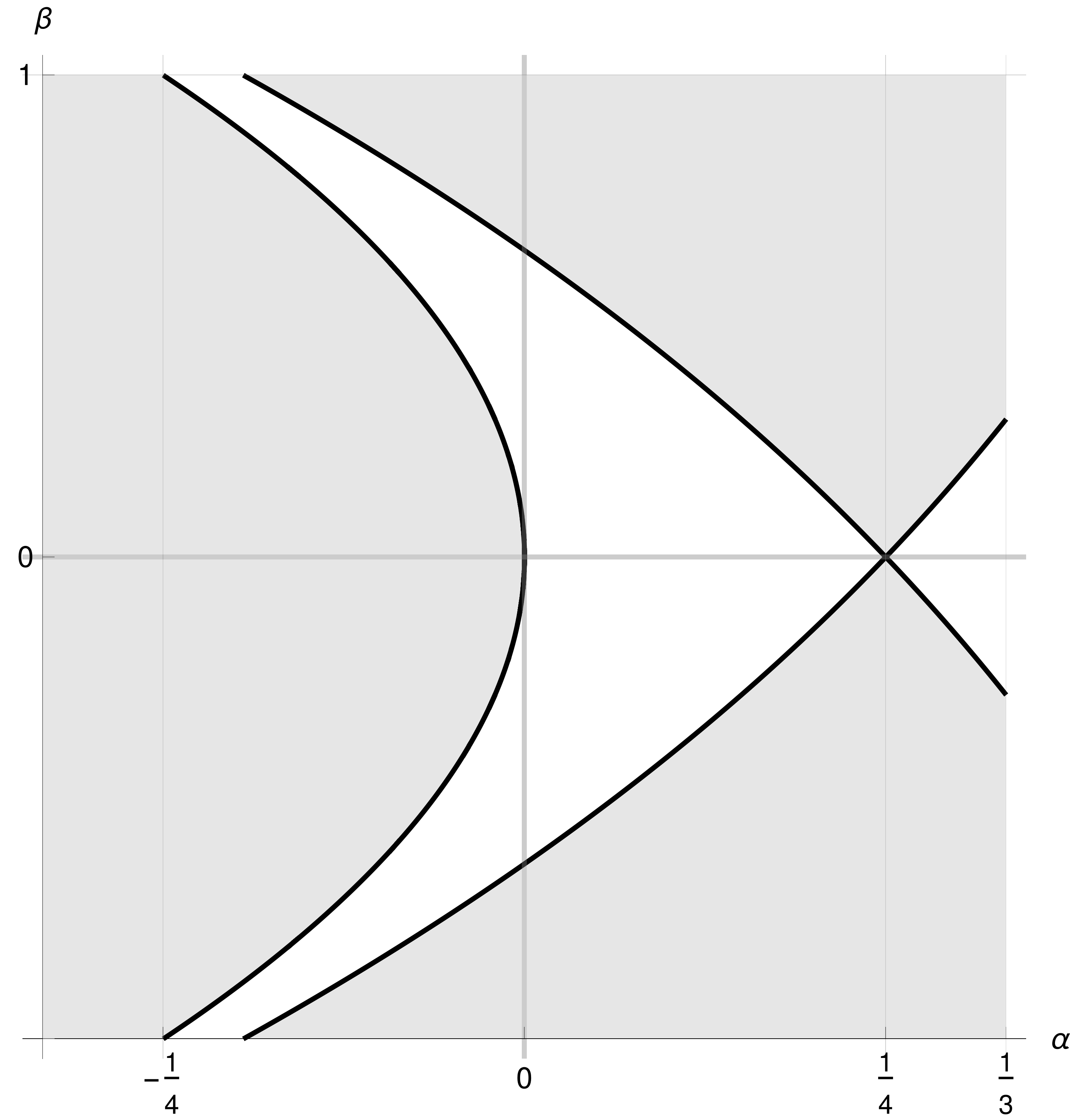}}
  \caption{Stable and unstable regions in the $(\alpha, \beta)$ plane ---  motion of pendulum pivot given by a triangular wave. Stable regions are shown in white, unstable regions are shown in grey.}
  \label{fig:triangular-wave-plots}
\end{figure}

In all cases discussed above the equations lead to explicit formulae describing curves in the parameter space $(\alpha, \beta) \in \R^2$, and it is straightforward to plot these curves and the corresponding stable/unstable regions in~$\R^2$, see Figure~\ref{fig:triangular-wave-plots}. Some important quantitative characteristics of the stability curves are easy to find from the obtained analytical expressions. For example, it is straightforward to verify that the curves emanating from the points $\begin{bmatrix} \alpha & \beta \end{bmatrix} = \begin{bmatrix} 0 & 0 \end{bmatrix}$ and $\begin{bmatrix} \alpha & \beta \end{bmatrix} = \begin{bmatrix} \frac{1}{4} & 0 \end{bmatrix}$ get closer as $\alpha \to -\infty$, which means that the stability gap in the negative half-plane is getting narrower as $\alpha \to -\infty$. Similarly the intersection points with the vertical axis are easy to find as well.

Finally, let us remark that since we have the analytical formula for the monodromy matrix, the results can be, if necessary, easily generalised to the case of \emph{damped} pendulum via the classical transformation, see, for example, \cite{jordan.dw.smith.p:nonlinear}.

\section{Pivot motion given as an approximated rectangular wave}
\label{sec:pivot-motion-given-2}

The stability analysis for the pivot motion given by an approximated rectangular wave goes along the same lines as for the triangular wave.
\subsection{Rectangular wave -- governing equations}
\label{sec:rect-wave-govern}
The pivot motion is, in the dimensionless time, given by the $2\pi$-periodic function $\xi_n$ specified by the formula
\begin{equation}
  \label{eq:44}
  \xi_n(\dimless{t})
  =
  \begin{cases}
    nA\dimless{t}, & \dimless{t} \in \left[-\frac{1}{n}, \frac{1}{n}\right], \\
    A, & \dimless{t} \in \left(\frac{1}{n}, \pi - \frac{1}{n} \right), \\
    -nA \left(\dimless{t} - \pi\right) & \dimless{t} \in \left[\pi - \frac{1}{n} , \pi + \frac{1}{n} \right], \\
    -A, & \dimless{t} \in \left(\pi + \frac{1}{n}, 2 \pi - \frac{1}{n} \right],
  \end{cases}
\end{equation}
see Figure~\ref{fig:pivot-motion-b}. This time we consider a sequence of possible pivot motions; the sequence labeled by $n \in \N$ approximates for $n \to + \infty$ the exact rectangular wave. The first and second time derivatives of $\xi$ are---informally---given by the formulae
\begin{equation}
  \label{eq:23}
  \dd{\xi_n}{\dimless{t}}
  =
  \begin{cases}
    nA, & \dimless{t} \in \left[-\frac{1}{n}, \frac{1}{n}\right], \\
    0, & \dimless{t} \in \left(\frac{1}{n}, \pi - \frac{1}{n} \right), \\
    -nA & \dimless{t} \in \left[\pi - \frac{1}{n} , \pi + \frac{1}{n} \right], \\
    0, & \dimless{t} \in \left(\pi + \frac{1}{n}, 2 \pi - \frac{1}{n} \right].
  \end{cases}
\end{equation}
and
\begin{equation}
  \label{eq:24}
  \ddd{\xi_n}{\dimless{t}}
  =
  nA \diracdelta_{\dimless{t} + \frac{1}{n}}
  -
  nA \diracdelta_{\dimless{t} - \frac{1}{n}}
  -
  nA \diracdelta_{\dimless{t} - \left(\pi - \frac{1}{n} \right)}
  +
  nA \diracdelta_{\dimless{t} - \left(\pi + \frac{1}{n} \right)}
  .
\end{equation}
As in the triangular wave case we see that the acceleration is given as a sequence of impulses, a sketch of the acceleration is given in Figure~\ref{fig:pivot-acceleration-b}. Substituting~\eqref{eq:24} into~\eqref{eq:6} and linearising with respect to $\dimless{\theta}$ yields
\begin{equation}
  \label{eq:25}
    \ddd{\dimless{\theta}}{\dimless{t}}
  +
  \left(
    \alpha
    +
    n
    \beta
    \left(
      \diracdelta_{\dimless{t} + \frac{1}{n}}
      -
      \diracdelta_{\dimless{t} - \frac{1}{n}}
      -
      \diracdelta_{\dimless{t} - \left(\pi - \frac{1}{n} \right)}
      +
      \diracdelta_{\dimless{t} - \left(\pi + \frac{1}{n} \right)}
    \right)
  \right)
  \dimless{\theta}
  =
  0
  ,
\end{equation}
where the dimensionless parameter $\beta$ is now given by the formula $\beta =_{\bydefinition} \frac{A}{l}$. The negative parameter values again correspond to the upright position, while the positive parameter values correspond to the pendent position.

\subsection{Monodromy matrix}
\label{sec:monodromy-matrix-1}


The monodromy matrix $\tensorq{E}$ for the system of interest~\eqref{eq:25} is the fundamental matrix evaluated over the whole period  $\tensorq{\Phi}_{-\frac{1}{n} -\to \left( 2 \pi - \frac{1}{n} \right) -}$, and it is again given by the product of the fundamental matrices for the free regime and the fundamental matrices for the single jumps located at the corresponding time instants, see~\eqref{eq:transfer-matrix-free} and~\eqref{eq:transfer-matrix-jump}. We get
\begin{multline}
  \label{eq:28}
  \tensorq{E}
  =
  \tensorq{\Phi}_{-\frac{1}{n} - \to \left(2 \pi - \frac{1}{n}\right) -}
  =
  \underbrace{
    \begin{bmatrix}
    \cos \left(  \sqrt{\alpha} \left(\pi - \frac{2}{n}\right) \right) &   \frac{1}{\sqrt{\alpha}} \sin \left(  \sqrt{\alpha} \left(\pi - \frac{2}{n}\right) \right) \\
    -\sqrt{\alpha} \sin \left(  \sqrt{\alpha} \left(\pi - \frac{2}{n}\right) \right) & \cos \left(  \sqrt{\alpha} \left(\pi - \frac{2}{n}\right) \right)
  \end{bmatrix}
}_{
  \tensorq{\Phi}_{\left(\pi + \frac{1}{n}\right)+ \to \left(2\pi - \frac{1}{n}\right)-}
}
  \\
  \times
  \underbrace{
    \begin{bmatrix}
      1 & 0 \\
      n\beta & 1 
    \end{bmatrix}
  }_
  {
     \tensorq{\Phi}_{\left(\pi + \frac{1}{n}\right)- \to \left(\pi + \frac{1}{n}\right)+}
  }
  \underbrace{
    \begin{bmatrix}
      \cos \left(  \sqrt{\alpha} \frac{2}{n} \right) &   \frac{1}{\sqrt{\alpha}} \sin \left(  \sqrt{\alpha} \frac{2}{n} \right) \\
      -\sqrt{\alpha} \sin \left(  \sqrt{\alpha} \frac{2}{n} \right) & \cos \left(  \sqrt{\alpha} \frac{2}{n} \right)
    \end{bmatrix}
  }_
  {
    \tensorq{\Phi}_{\left(\pi - \frac{1}{n}\right)+ \to \left(\pi + \frac{1}{n}\right)-}
  }
  \underbrace{
    \begin{bmatrix}
      1 & 0 \\
      -n\beta & 1 
    \end{bmatrix}
  }_
  {
    \tensorq{\Phi}_{\left(\pi - \frac{1}{n}\right)- \to \left(\pi - \frac{1}{n}\right)+}
  }
  \\
  \times
  \underbrace{
    \begin{bmatrix}
      \cos \left(  \sqrt{\alpha} \left(\pi - \frac{2}{n}\right) \right) &   \frac{1}{\sqrt{\alpha}} \sin \left(  \sqrt{\alpha} \left(\pi - \frac{2}{n}\right) \right) \\
      -\sqrt{\alpha} \sin \left(  \sqrt{\alpha} \left(\pi - \frac{2}{n}\right) \right) & \cos \left(  \sqrt{\alpha} \left(\pi - \frac{2}{n}\right) \right)
    \end{bmatrix}
  }
  _
  {
    \tensorq{\Phi}_{\frac{1}{n}+ \to \left(\pi - \frac{1}{n}\right)-}
  }
  \\
  \times
  \underbrace{
    \begin{bmatrix}
      1 & 0 \\
      -n\beta & 1 
    \end{bmatrix}
  }_
  {
    \tensorq{\Phi}_{\frac{1}{n}- \to \frac{1}{n}+}
  }
  \underbrace{
    \begin{bmatrix}
      \cos \left(  \sqrt{\alpha} \frac{2}{n} \right) &   \frac{1}{\sqrt{\alpha}} \sin \left(  \sqrt{\alpha} \frac{2}{n} \right) \\
      -\sqrt{\alpha} \sin \left(  \sqrt{\alpha} \frac{2}{n} \right) & \cos \left(  \sqrt{\alpha} \frac{2}{n} \right)
    \end{bmatrix}
  }_
  {
    \tensorq{\Phi}_{-\frac{1}{n}+ \to \frac{1}{n}-}
  }
  \underbrace{
    \begin{bmatrix}
      1 & 0 \\
      n\beta & 1 
    \end{bmatrix}
  }_
  {
     \tensorq{\Phi}_{-\frac{1}{n}- \to -\frac{1}{n}+}
  }
  ,
\end{multline}
which upon the application of trigonometric identities yields
\begin{multline}
  \label{eq:29}
  \Tr
  \tensorq{\Phi}_{-\frac{1}{n} - \to \left(2 \pi - \frac{1}{n}\right)-}
\\
=
\frac{4 \beta ^2 n^2 \left(4 \alpha +\beta ^2 n^2\right) \sin ^2\left(\pi \sqrt{\alpha} - \frac{2\sqrt{\alpha}}{n}\right)
  +
  \beta ^4 n^4 \cos \left(2 \pi \sqrt{\alpha} - \frac{8 \sqrt{\alpha}}{n}\right)
  -
  2 \beta ^2 n^2 \left(4 \alpha +\beta ^2 n^2\right) \cos \left(\frac{4 \sqrt{\alpha }}{n}\right)
  +
  \cos \left(2 \pi  \sqrt{\alpha }\right) \left(4 \alpha +\beta ^2 n^2\right)^2}{8 \alpha ^2}
\\
=
n^2
\left(
  \frac{4 \beta^2 \sin^2 \left( \pi \sqrt{\alpha} \right)}{\alpha}
  -
  \frac{
    4 \beta^2
    \left(
      1
      +
      2 \beta^2
    \right)
    \sin \left(2 \pi \sqrt{\alpha} \right)
  }
  {
    \sqrt{\alpha}
  }
  \frac{1}{n}
  +
  \left[
    \left(\frac{56 \beta ^4}{3}+8 \beta ^2+2\right) \cos \left(2 \pi  \sqrt{\alpha }\right)-\frac{8}{3} \beta ^2 \left(\beta ^2-3\right)
  \right]
  \frac{1}{n^2}
  +
  \bigo{\frac{1}{n^3}}
\right)
.
\end{multline}
Note that the limit $n \to + \infty$ is not finite, hence we can not directly analyse the exact rectangular wave, we must stay on the level of approximated rectangular wave.

\subsection{Search for periodic solutions}
\label{sec:search-peri-solut-1}
Following again the standard analysis, see~\cite{jordan.dw.smith.p:nonlinear}, we know that the system~\eqref{eq:8} possesses a periodic solution provided that $\Tr \tensorq{\Phi}_{-\frac{1}{n} - \to \left(2 \pi - \frac{1}{n}\right)-} = \pm 2$, and that the implicitly defined curves $\Tr \tensorq{\Phi}_{-\frac{1}{n} - \to \left(2 \pi - \frac{1}{n}\right)-} = \pm 2$ separate the stable/unstable regions in the parameter space.

\subsubsection{Solution of equation $\Tr \tensorq{\Phi}_{-\frac{1}{n} - \to \left(2 \pi - \frac{1}{n}\right)-} = 2$}
\label{sec:solution-equation-tr-2}
\begin{figure}[h]
  \centering
  \subfloat[$n=4$ \label{fig:rectangular-wave-plots-a}]{\includegraphics[width=0.32\textwidth]{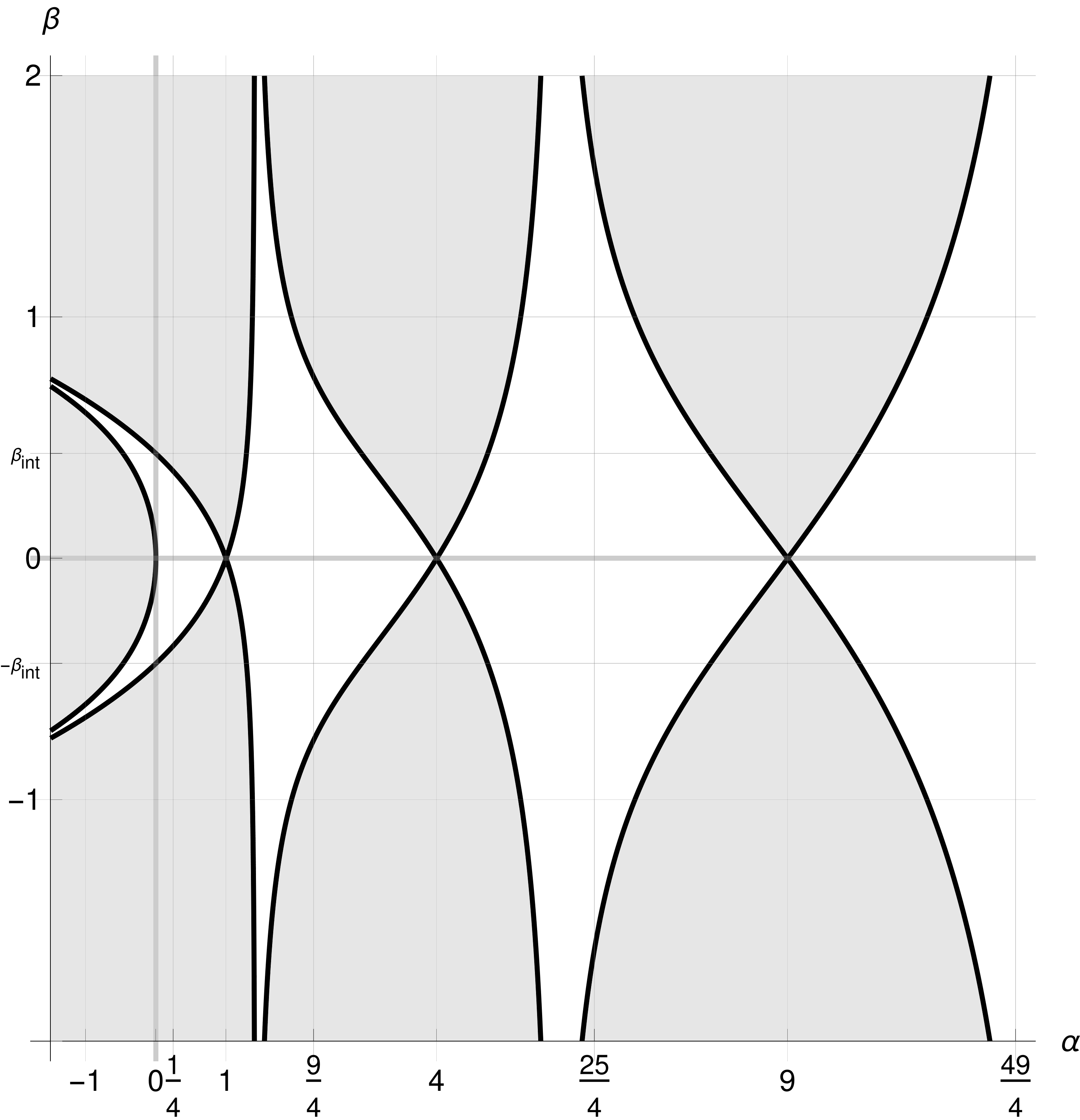}}
  \quad
  \subfloat[$n=10$ \label{fig:rectangular-wave-plots-b}]{\includegraphics[width=0.32\textwidth]{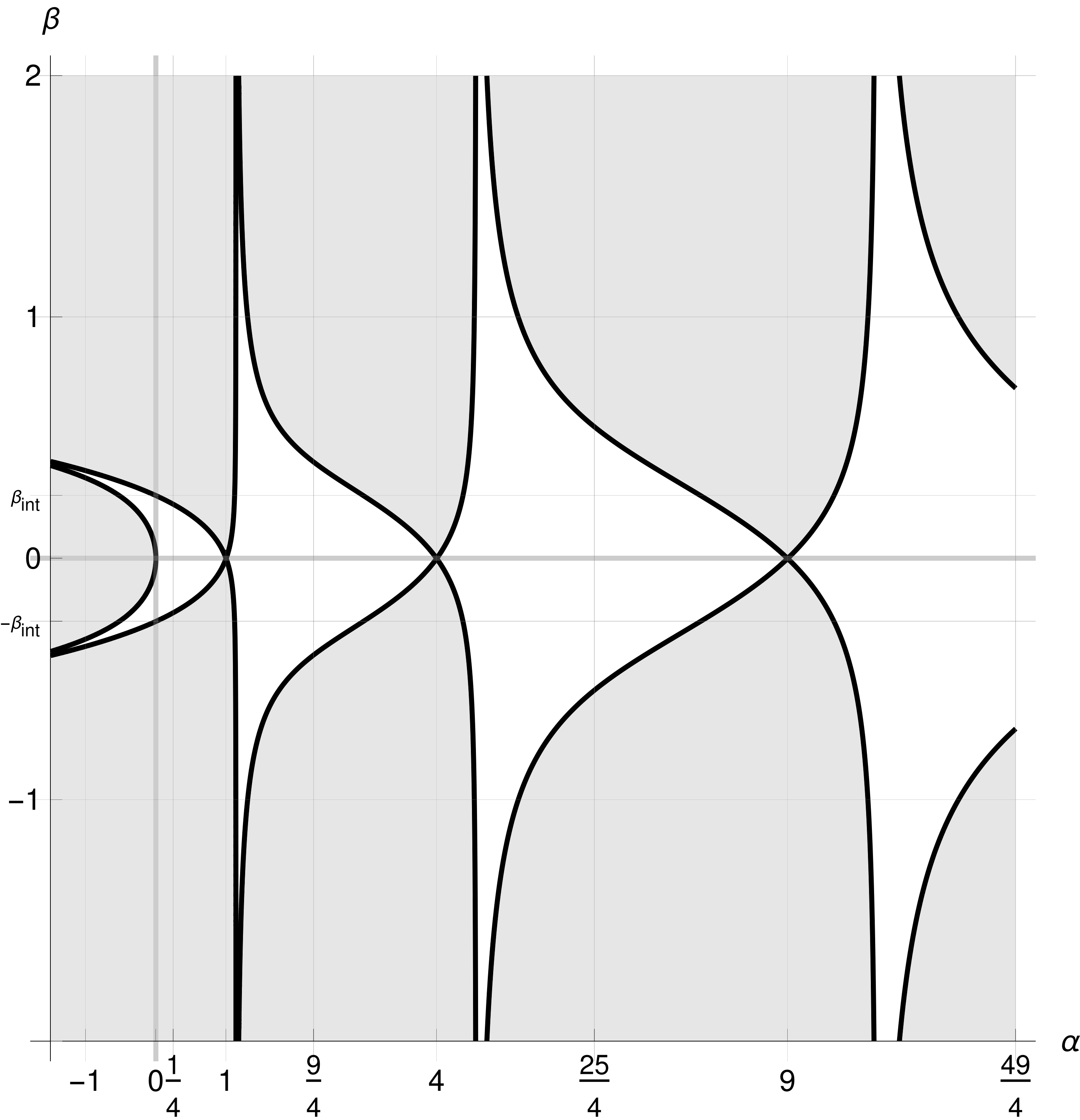}}
  \\
  \subfloat[$n=20$ \label{fig:rectangular-wave-plots-c}]{\includegraphics[width=0.32\textwidth]{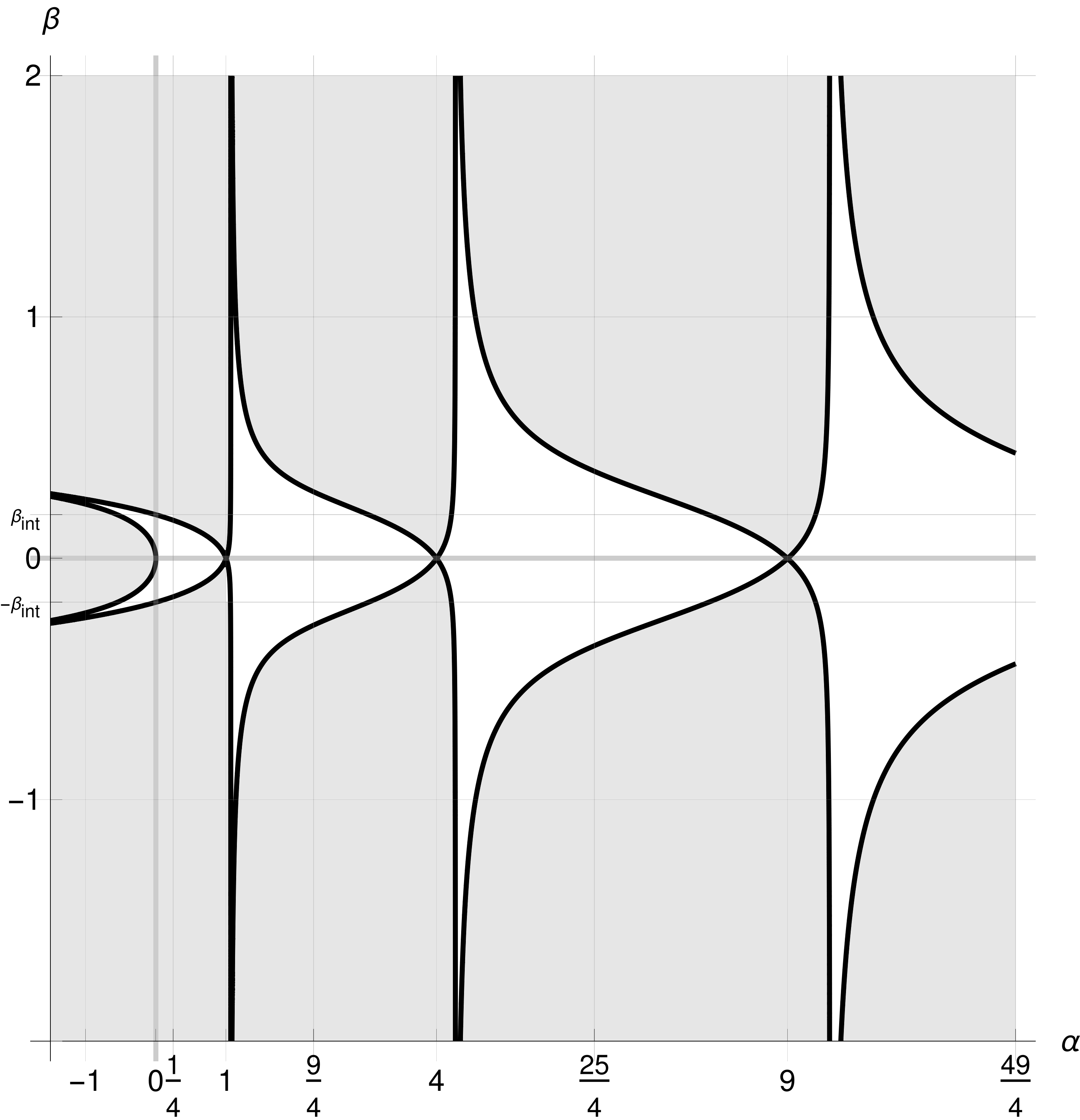}}
  \quad
  \subfloat[$n=100$ \label{fig:rectangular-wave-plots-d}]{\includegraphics[width=0.32\textwidth]{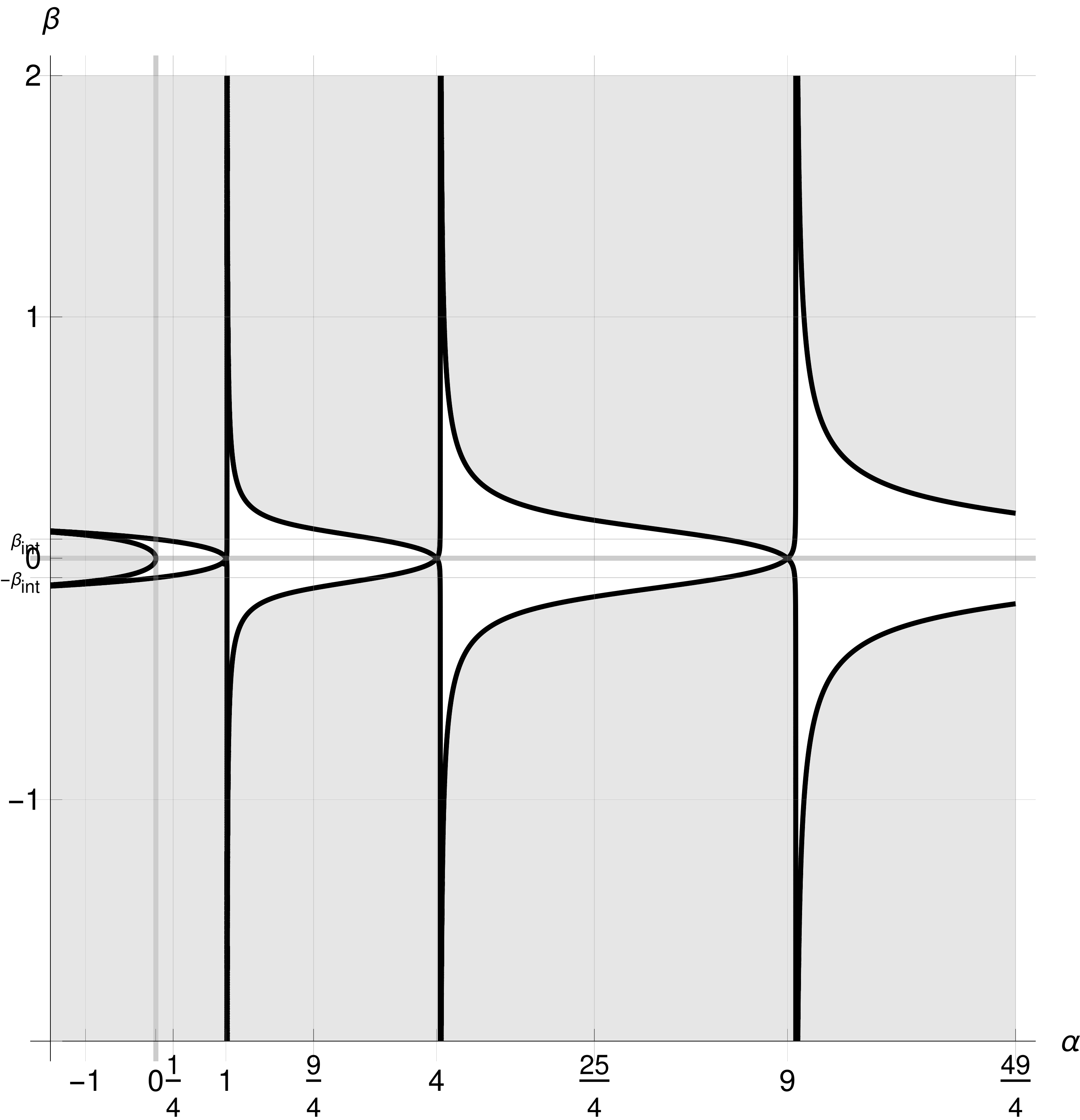}}
  \caption{Stable and unstable regions in the $(\alpha, \beta)$ plane ---  motion of pendulum pivot given by an approximation of a rectangular wave. Stable regions are shown in white, unstable regions are shown in grey.}
  \label{fig:rectangular-wave-plots}
\end{figure}

Solutions to this equation can be without any difficulties found numerically, and the correspondig Ince--Strutt diagram is straightforward to plot, see Figure~\ref{fig:rectangular-wave-plots}. Some quantitative characteristics are however given by explicit formulae. For example, the intersections of the solution curves with the horizontal axis $\beta = 0$ are located at points where $\left. \left( \Tr \tensorq{\Phi}_{-\frac{1}{n} - \to \left(2 \pi - \frac{1}{n}\right)-} \right) \right|_{\beta = 0} = 2$, which leads to the equation
$
2 \cos \left(2 \pi \sqrt{\alpha} \right) = 2
$
with solutions
\begin{equation}
  \label{eq:34}
  \alpha = k^2,
\end{equation}
where $k \in \N$. Nontrivial intersection points with the vertical axis $\alpha=0$ are also straightforward to find, $\beta_{\text{int}}=\pm \sqrt{\frac{2}{n \pi -2}}$. 

\subsubsection{Solution of equation $\Tr \tensorq{\Phi}_{-\frac{1}{n} - \to \left(2 \pi - \frac{1}{n}\right)-} = -2$}
\label{sec:solution-equation-tr-3}
If $\beta=0$, then we use~\eqref{eq:29}, and we see that the equation $\Tr \tensorq{\Phi}_{-\frac{1}{n} - \to \left(2 \pi - \frac{1}{n}\right)-} = -2$ reduces to
$
2 \cos \left(2 \pi \sqrt{\alpha} \right) = -2
$.
This equation has solutions
\begin{equation}
  \label{eq:32}
  \alpha = \left(k + \frac{1}{2}\right)^2,
\end{equation}
where $k \in \N$. On the other hand, if $\beta \not = 0$ and $n$ is sufficiently high, then the equation  $\Tr \tensorq{\Phi}_{-\frac{1}{n} - \to \left(2 \pi - \frac{1}{n}\right)-} = -2$ has \emph{no solutions}. The observation that the equation $\Tr \tensorq{\Phi}_{-\frac{1}{n} - \to \left(2 \pi - \frac{1}{n}\right)-} = -2$ has only trivial solutions with $\beta = 0$ constitutes the main qualitative feature of the system response to the forcing by an approximated rectangular wave. The classical subharmonic resonances at $\alpha = \frac{1}{4}$ and other frequencies given by the formula~\eqref{eq:32} \emph{are not present}. 

The fact that for sufficiently high $n$ there are no additional solutions in the neighbourhood of the point
$
\begin{bmatrix}
\alpha & \beta  
\end{bmatrix}
=
\begin{bmatrix}
  \frac{1}{4} & 0
\end{bmatrix}
$
can be shown by the following argument, which is easy to extend to other $\alpha$ parameters in the form~\eqref{eq:32}. The equation $\Tr \tensorq{\Phi}_{-\frac{1}{n} - \to \left(2 \pi - \frac{1}{n}\right)-} = -2$ 
can be rewritten as
\begin{multline}
  \label{eq:37}
  \beta^4 n^4
  \underbrace{
    \left[
      4
      \sin ^2\left(\pi \sqrt{\alpha} - \frac{2\sqrt{\alpha}}{n}\right)
      +
      \cos \left(2 \pi \sqrt{\alpha} - \frac{8 \sqrt{\alpha}}{n}\right)
      -
      2
      \cos \left(\frac{4 \sqrt{\alpha }}{n}\right)
      +
      \cos \left(2 \pi \sqrt{\alpha}\right)
    \right]
  }_{A}
  \\
  +
  4 \alpha \beta^2n^2
  \underbrace{
    \left[
      4
      \sin ^2\left(\pi \sqrt{\alpha} - \frac{2\sqrt{\alpha}}{n}\right)
      -
      2
      \cos \left(\frac{4 \sqrt{\alpha }}{n}\right)
      +
      2\cos \left(2 \pi \sqrt{\alpha}\right)
    \right]
  }_{B}
  =
  - 16 \alpha^2
  \left[
    1 + \cos \left(2 \pi \sqrt{\alpha}\right)
  \right]
  .
\end{multline}
Clearly, the pair $\alpha = \frac{1}{4}$ and $\beta = 0$ is a solution to~\eqref{eq:37} for any $n \in \N$. On the other hand, if $\alpha \approx \frac{1}{4}$, meaning that if $\alpha$ is close to~$\frac{1}{4}$ but not equal to~$\frac{1}{4}$, then the right-hand side of~\eqref{eq:37} is negative, while the terms $A$ and $B$ are for sufficiently high $n$ nonnegative. Consequently, the equation has no solution in the neighbourhood of $\frac{1}{4}$ except of the trivial solution $\alpha = \frac{1}{4}$ and~$\beta = 0$. The nonnegativity of $A$ follows from the manipulation
\begin{multline}
  A
  =
  4
  \sin ^2\left(\pi \sqrt{\alpha} - \frac{2\sqrt{\alpha}}{n}\right)
  +
  \cos \left(2 \pi \sqrt{\alpha} - \frac{8 \sqrt{\alpha}}{n}\right)
  -
  2
  \cos \left(\frac{4 \sqrt{\alpha }}{n}\right)
  +
  \cos \left(2 \pi \sqrt{\alpha}\right)
  \\
  =
  2
  \left[
    1
    -
    \cos \left(\frac{4 \sqrt{\alpha}}{n}\right)
  \right]
  \left[
    1
    -
    \cos \left(2 \pi \sqrt{\alpha} - \frac{4 \sqrt{\alpha}}{n} \right)
  \right]
  \geq 0
  ,
\end{multline}
and the nonnegativity holds for arbitrary $n$ and arbitrary positive $\alpha$. Regarding the nonnegativity of the term $B$, we see that the term can be rewritten as
\begin{equation}
  \label{eq:38}
  B
  =
  4
  \sin ^2\left(\pi \sqrt{\alpha} - \frac{2\sqrt{\alpha}}{n}\right)
  -
  2
  \cos \left(\frac{4 \sqrt{\alpha }}{n}\right)
  +
  2
  \cos \left(2 \pi \sqrt{\alpha}\right)
  =
  -
  8
  \sin
  \left(
    \frac{4\sqrt{\alpha}}{n}
  \right)
  \cos
  \left(
    2\pi \sqrt{\alpha}
  \right)
  \sin \left(\pi \sqrt{\alpha} - \frac{2\sqrt{\alpha}}{n}\right)
  \geq 0,
\end{equation}
where we need to consider sufficiently high $n$ and $\alpha$ close to $\frac{1}{4}$.

\section{Pivot motion given as a cosine wave -- classical case}
\label{sec:pivot-motion-given-3}
\begin{figure}[h]
  \centering
  \subfloat[Global view. \label{fig:cosine-wave-plots-a}]{\includegraphics[height=0.32\textwidth]{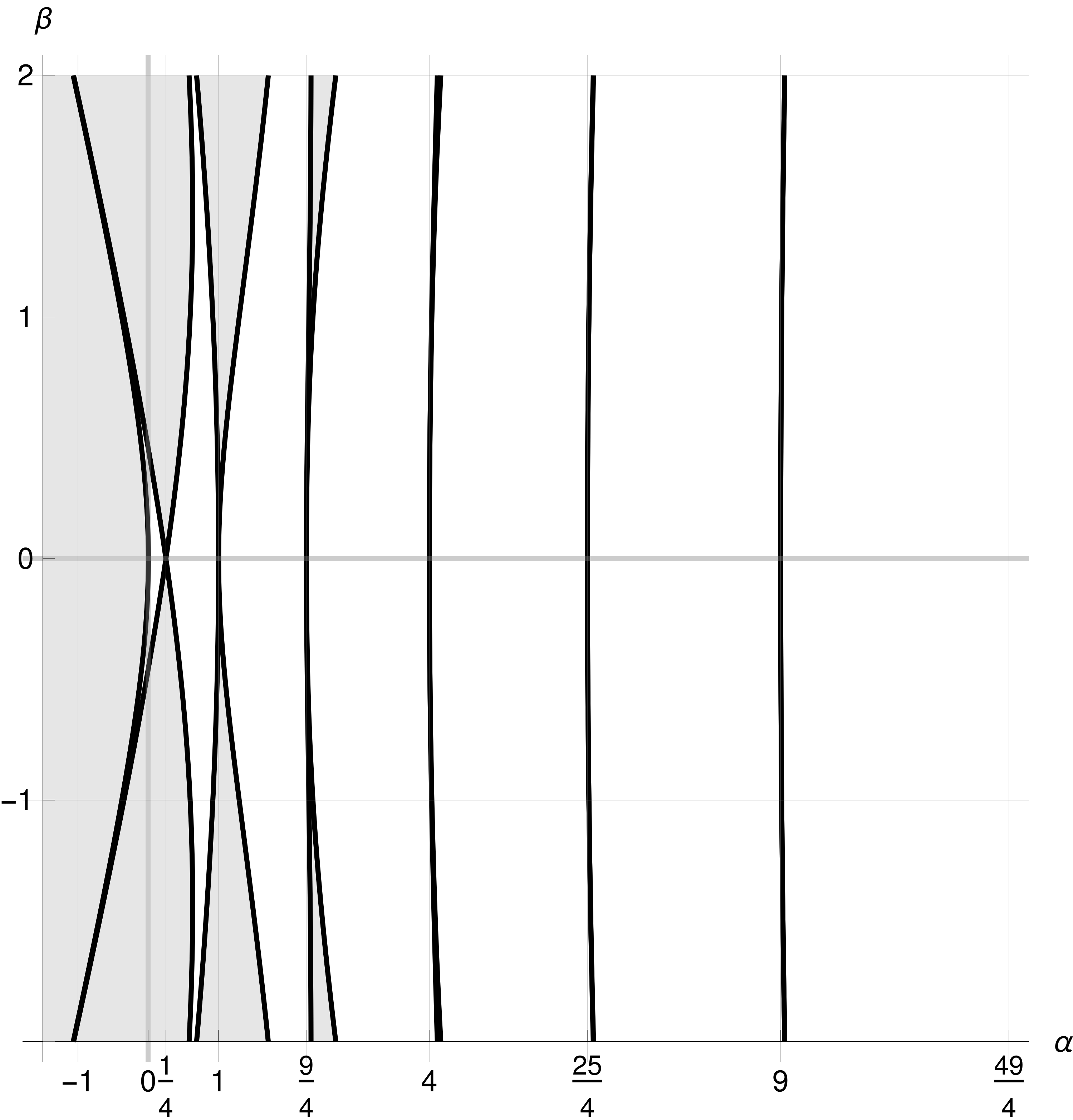}}
  \quad
  \subfloat[Detailed view. \label{fig:cosine-wave-plots-b}]{\includegraphics[height=0.32\textwidth]{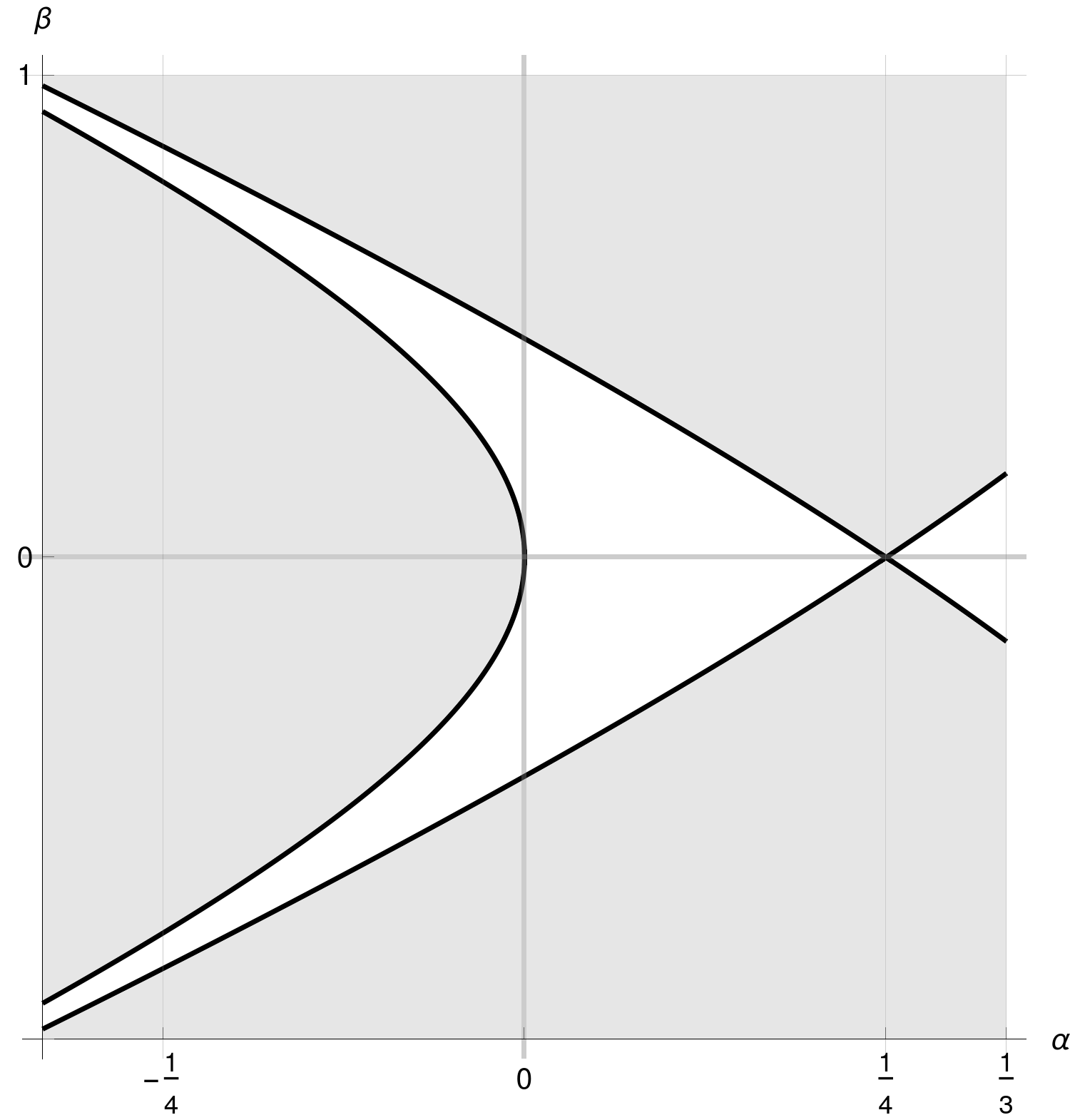}}
  \quad
  \subfloat[Global view -- wider range of $\beta$. \label{fig:cosine-wave-plots-c}]{\includegraphics[height=0.32\textwidth]{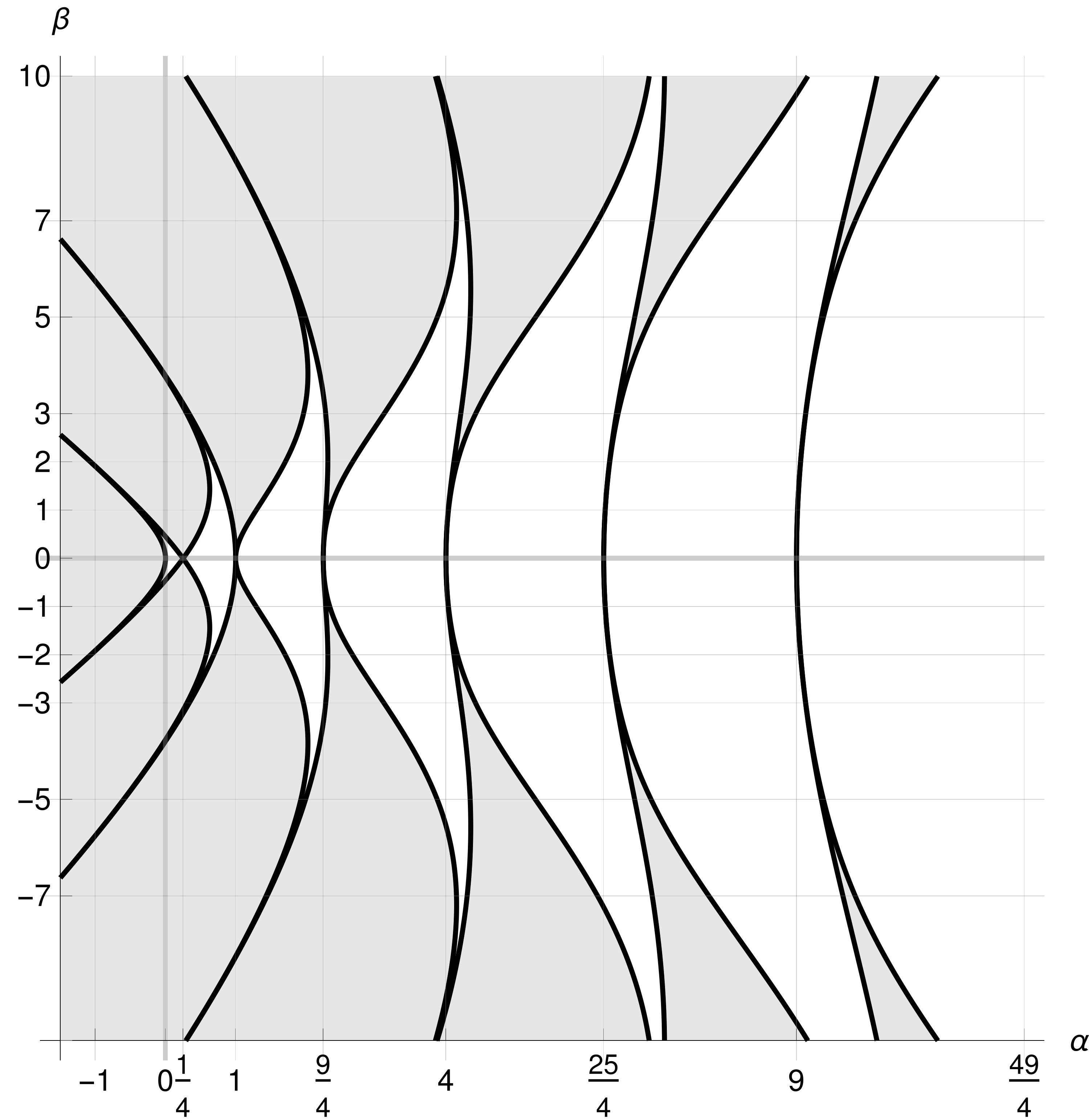}}
  \caption{Classical Ince--Strutt diagram. Stable and unstable regions in the $(\alpha, \beta)$ plane -- motion of pendulum pivot given by a cosine wave. Stable regions are shown in white, unstable regions are shown in grey.}
  \label{fig:cosine-wave-plots}
\end{figure}

For the sake of completeness and subsequent discussion we also show the classical Ince--Strutt diagram for the pivot motion given by a cosine wave of angular frequency $\Omega$ and amplitude $A$, $\xi = A \cos (\Omega t)$. The cosine wave has the same amplitude and period as the triangular wave and the approximated rectangular wave investigated in the previous sections. The corresponding linearised dimensionless governing equation is the Mathieu equation~\eqref{eq:39}, where $\alpha=_{\bydefinition}\frac{\omega_0^2}{\Omega^2}$ and  $\beta =_{\bydefinition} -\frac{A}{l}$, and standard stability analysis leads to the well-known Ince--Strutt stability diagram shown in Figure~\ref{fig:cosine-wave-plots}. Note that the diagram is nowadays easy to produce, since many software packages, for example~\texttt{Wolfram Mathematica}, have built-in routines for special functions such as Mathieu functions. In particular, the diagram shown in Figure~\ref{fig:cosine-wave-plots} has been produced in \texttt{Wolfram Mathematica} using the built-in routines \texttt{MathieuCharacteristicA} and \texttt{MathieuCharacteristicB}. However, the implementation of routines for special functions such as Mathieu functions is far from being trivial, see, for example, \cite{alhargan.fa:complete} for a related discussion, and it usually heavily relies on fine analytical results and approximations, see especially the pioneering work by \cite{ince.el:researches} and numerous monographs on Mathieu functions such as~\cite{mclachlan.nw:theory} or~\cite{magnus.wws:hills}.

\section{Discussion and conclusion}
\label{sec:disc-concl}

Stability diagrams for the mathematical pendulum with the pivot motion described by the triangular wave and by the approximated rectangular wave are straightforward to obtain. In particular, the curves separating the stable and unstable regions in the stability diagram for the triangular wave are described by the explicit formula given in terms of elementary functions. This contrasts with the classical analysis for the pivot motion described by a cosine wave, that heavily relies on involved technical manipulations and various properties of special functions.

Interestingly, basic qualitative features of the stability diagram for the~\emph{triangular wave} are the same as in the classical case of a cosine wave. Indeed, both stability diagrams have a stability region in the negative half-plane and close to the origin, compare Figure~\ref{fig:triangular-wave-plots} and Figure~\ref{fig:cosine-wave-plots}. Consequently, both pivot motions can stabilise the pendulum in the upright position. Furthermore, the classical subharmonic resonance region close to $\alpha = \frac{1}{4}$ is present in both stability diagrams, compare again Figure~\ref{fig:triangular-wave-plots} and Figure~\ref{fig:cosine-wave-plots}. Since the qualitative behaviour of the two systems is in this regard the same, and since the rigorous analysis for the triangular wave is much simpler and complete than in the classical case of the cosine wave, it seems that the triangular wave pivot motion is an ideal candidate for the instructive presentation of the parametric resonance phenomenon in the case of pendulum with a moving pivot.

The approximated \emph{rectangular wave} also leads to a stability region in the negative half-plane and close to the origin, compare Figure~\ref{fig:rectangular-wave-plots} and Figure~\ref{fig:cosine-wave-plots}, hence even this pivot motion can stabilise the pendulum in the upright position. However, the counterpart of the classical subharmonic resonance close to $\alpha = \frac{1}{4}$ is in this case \emph{not present}, which challenges the classical ``rule of thumb'' regarding the destabilisation inducing frequencies via the parametric resonance mechanism.

The findings of our analysis might be transferred to more complicated settings wherein the parametric resonance/excitation plays a role---see~\cite{champneys.a:dynamics} for a thorough review---and wherein simple and instructive analytical formulae are of interest.


\bibliography{vit-prusa}
\end{document}